# Transmission Expansion Planning with Seasonal Network Optimization

Xingpeng Li
Department of Electrical and Computer Engineering
University of Houston
Email: xli82@uh.edu

Qianxue Xia
School of Electrical and Computer Engineering
Georgia Institute of Technology
Email: qxia31@gatech.edu

*Abstract*— Transmission expansion planning (TEP) is critical for the power grid to meet fast growing demand in the future. Traditional TEP model does not utilize the flexibility in the transmission network that is considered as static assets. However, as the load profile may have different seasonal patterns, the optimal network configuration could be very different for different seasons in the planning horizon. Therefore, this paper proposes to incorporate seasonal network optimization (SNO) into the traditional TEP model. SNO dynamically optimizes the network for each season of each planning epoch. Two TEP-SNO models are proposed to investigate the benefits of optimizing the status of (i) existing branches, and (ii) existing and new branches, respectively. Numerical simulations demonstrate the effectiveness of the proposed TEP-SNO models. It is shown that SNO can improve system operational efficiency, defer investment of new transmission elements, and reduce the total cost.

## NOMENCLATURE

*Sets*:
- $E$   Epochs in the planning horizon.
- $G$   Generators.
- $G(n)$   Generators at bus $n$.
- $K$   Branches.
- $K(n-)$   Branches of which bus $n$ is the from-bus.
- $K(n+)$   Branches of which bus $n$ is the to-bus.
- $N$   Buses.
- $J$   Candidate new branches.
- $J(n-)$   Candidate new branches of which bus $n$ is from-bus.
- $J(n+)$   Candidate new branches of which bus $n$ is the to-bus.
- $S$   Seasons in a year.
- $T$   Hours in a day.

*Indices*:
- $e$   Epoch $e$, an element of set $E$.
- $g$   Generator $g$, an element of set $G$.
- $n$   Bus $n$, an element of set $N$.
- $j$   Candidate new branch $j$, an element of set $J$.
- $k$   Existing branch $k$, an element of set $K$.
- $s$   Season $s$, an element of set $S$.
- $t$   Hour $t$, an element of set $T$.

*Parameters*:
- $A_D$   Annual load increase in percent.
- $A_M$   Annual maintenance cost in percent for new branches.
- $C_g$   Variable production cost for unit $g$.
- $C_j$   Capital cost for constructing new branch $j$.
- $M$   A very big constant number.
- $N_E$   Number of epochs in the planning horizon.
- $N_{YE}$   Number of years in a planning epoch.
- $P_{g,max}$   Maximum output of unit $g$.
- $P_{g,min}$   Minimum output of unit $g$.
- $D_{ntse}$   Load at bus $n$ for hour $t$ in a typical day of season $s$ in epoch $e$.
- $RateA_k$   Long-term limit of existing branch $k$.
- $RateA_j$   Long-term limit of candidate new branch $j$.
- $TC$   Total cost for transmission investment and generation in the entire planning horizon.
- $TC_G$   Total generation cost in the planning horizon.
- $TC_I$   Total investment cost in the planning horizon.
- $x_k$   Reactance of existing branch $k$.
- $x_j$   Reactance of candidate new branch $j$.
- $n_{(k-)}$   From-bus of existing branch $k$.
- $n_{(k+)}$   To-bus of existing branch $k$.
- $n_{(j-)}$   From-bus of candidate new branch $j$.
- $n_{(j+)}$   To-bus of candidate new branch $j$.

*Variables*:
- $p_{gtse}$   Output power of unit $g$ in hour $t$ in a typical day of season $s$ in epoch $e$.
- $\theta_{ntse}$   Phase angle of bus $n$ in hour $t$ in a typical day of season $s$ in epoch $e$.
- $p_{ktse}$   Power flow on existing line $k$ in hour $t$ in a typical day of season $s$ in epoch $e$.
- $p_{jtse}$   Power flow on new line $j$ in hour $t$ in a typical day of season $s$ in epoch $e$.
- $u_{je}$   Availability of candidate new line $j$ in epoch $e$.
- $v_{je}$   Construction flag of candidate new line $j$ in epoch $e$.
- $z_{kse}$   Status of existing line $k$ in a typical day of season $s$ in epoch $e$.
- $z_{jse}$   Status of candidate new line $j$ in a typical day of season $s$ in epoch $e$.

## I. INTRODUCTION

Electricity consumption increases overtime all over the world. The electrical network that is used to deliver electric power produced at remote power plants to load pockets has limited capacity, which may prevent cheap units from producing at their capacities. Moreover, network congestion would become much more severe when load reaches a significantly higher level in the next few years. Thus, long-term transmission expansion planning (TEP) is required to create a cost-effective plan for upgrading existing transmission network to efficiently meet the load growth in the future.

Since the capital cost of constructing a new transmission component (line or transformer) is very high and the constructed new line will serve for a great number of years, it is critical for developing an effective procedure to make such investment decisions. Typically, a co-optimization model that incorporates both transmission investment plan and generation scheduling is formulated and solved to achieve least cost solutions for TEP problems [1]-[2]. A Benders decomposition approach is proposed in [3] to solve transmission network design problems; this work developed a strategy to minimize the value of the big-M. Reference [4] proposes a TEP model that considers reactive power and power losses. A heuristic is developed in [5] to reduce the combinatorial search space of the TEP problem such that the TEP problem can be solved with a low computational time. A multi-objective TEP framework that can provide Pareto optimal solutions is proposed in [6]. Reference [7] proposes a two-stage adaptive robust TEP model that considers netload uncertainty. The placement of battery energy storage systems is co-optimized with TEP in [8]; it shows that inclusion of battery energy storage systems can defer the construction of some new lines in some scenarios.

Though the flexibility in the transmission network has not been fully utilized in practical power system operations, prior efforts have shown that dynamically optimizing the network configuration can achieve cost savings [9]-[13], congestion relief [14]-[15], and reliability enhancement [16]-[18]. However, these studies are for short-term operations and very few efforts demonstrate the benefits of optimizing the transmission network in the long-term TEP problems.

Traditionally, existing transmission elements are assumed to be connected in the network all the time during the entire planning horizon. However, the system load profile may not share the same patterns in different seasons; as a result, the optimal network configuration could be different for different seasons. Moreover, due to annual load growth, the optimal network topology could be different for the same type of season in different years. Therefore, seasonal network optimization (SNO) may determine more efficient network configurations for different seasons in different planning epochs and incorporating SNO into TEP may reduce the total cost.

In this paper, a traditional TEP model is first established, followed by the transmission expansion planning with seasonal network optimization (TEP-SNO) model. Two types of TEP-SNO models are proposed to investigate the benefits of optimizing the status of branches. TEP-SNO-T1 only adjusts the status of branches that already exist in the network before the planning while TEP-SNO-T2 also optimizes the status of new lines after their construction is complete. Numerical simulations demonstrate that including SNO in TEP problems can eliminate unnecessary transmission investment, defer the construction of new lines and improve network efficiency; thus, it can significantly decrease the total expected cost.

The rest of this paper is organized as follows: Section II presents the proposed methodology and model; Section III analyzes the results of numerical simulations; and Section IV concludes the paper.

## II. METHODOLOGY AND MODEL

Long-term transmission expansion planning is critical for meeting the growing demand and maintaining the system reliability for the future electricity grid. This section will introduce and explain the detailed formulations for the traditional TEP model and the proposed two enhanced TEP-SNO models.

### A. TEP

The objective function of TEP presented in (1) is to minimize the total expected cost $TC$.

$$min\ TC = TC_G + TC_I \quad (1)$$

where $TC_G$ denotes the total generation cost in the planning horizon; and $TC_I$ denotes the total capital cost for building new lines. The equations for calculating $TC_G$ and $TC_I$ are presented in (2) and (3) respectively.

$$TC_G = N_{YE} \frac{365}{4} \sum_g \sum_t \sum_s \sum_e p_{gtse} C_g \quad (2)$$

$$TC_I = \sum_j \sum_e v_{je} C_j (1 + (N_E - e + 1) N_{YE} A_M) \quad (3)$$

Each season is represented by a typical day; in other words, the annual load profile is simplified and represented by four typical days that represent Spring, Summer, Fall and Winter respectively. This explains the existence of term $\frac{365}{4}$ in (2). The annual percentage maintenance cost for each candidate new branch is assumed to be the same, 4%. This annual maintenance cost may include transmission degradation cost, crew labor cost and other related direct/indirect cost.

Nodal power balance constraints are enforced in (4) for each hourly interval in each season of each planning epoch. Equation (5) calculate the load for each epoch considering annual load growth. In this paper, it is assumed that the annual load growth rate $A_D$ is 2% and the annual load profile remains the same for the entire epoch period. The beginning-year load profile is used for all the years in the same epoch. Generator production limits are modeled in (6). To avoid modeling unit commitment that would significantly increase the computational burden, the minimum generation is set to be zero for each unit in this work. Transmission element capacity limit is enforced in (7). Equation (8) defines the relationship between bus phase angles and branch power flows.

$$\sum_{k \in K(n+)} p_{ktse} - \sum_{k \in K(n-)} p_{ktse} + \sum_{j \in J(n+)} p_{jtse} -$$
$$\sum_{j \in J(n-)} p_{jtse} + \sum_{g \in G(n)} p_{gtse} = D_{ntse} \quad \forall n, t, s, e \quad (4)$$

$$D_{ntse} = D_{nts1} (1 + A_D)^{(e-1)N_{YE}} \quad \forall n, t, s, e \quad (5)$$

$$P_{g,min} \le p_{gtse} \le P_{g,max} \quad \forall g, t, s, e \quad (6)$$

$$-RateA_k \leq p_{ktse} \leq RateA_k \quad \forall k,t,s,e \quad (7)$$

$$p_{ktse} = (\theta_{n_{(k-)}tse} - \theta_{n_{(k+)}tse})/x_k \quad \forall k,t,s,e \quad (8)$$

To incorporate new line investment in the TEP model, two binary variables are defined: $v_{je}$ and $u_{je}$. $v_{je}$ denotes the construction flag of candidate new line $j$ in epoch $e$. When $v_{je}$ is 1, it means the new line with an index of $j$ will be constructed and finished at the beginning of epoch $e$; otherwise, it is 0. $u_{je}$ denotes the availability of candidate new line $j$ in epoch $e$. When $u_{je}$ is 1, it denotes that the new line with an index of $j$ is available for use in the system, which indicates the new line $j$ has been completely built at or before the beginning of epoch $e$. The power flow on candidate new line $j$ can be determine by the *Big-M* method, as shown in (9) and (10). Constraint (11) models the capacity limits for candidate new lines. Constraints (9)-(11) together can ensure that (i) after its construction, a new line will respect the same physical restrictions as existing lines, and (ii) a new line before its construction will not affect the network flow. Constraints (12)-(14) define the relationship between $v_{je}$ and $u_{je}$.

$$p_{jtse} - \frac{\theta_{n_{(j-)}tse} - \theta_{n_{(j+)}tse}}{x_j} \leq (1 - u_{je})M \quad \forall k,t,s,e \quad (9)$$

$$p_{jtse} - \frac{\theta_{n_{(j-)}tse} - \theta_{n_{(j+)}tse}}{x_j} \geq -(1 - u_{je})M \quad \forall k,t,s,e \quad (10)$$

$$-u_{je}RateA_j \leq p_{jtse} \leq u_{je}RateA_j \quad \forall j,t,s,e \quad (11)$$

$$\sum_{z \leq e} v_{jz} \leq u_{je} \quad \forall j,e \quad (12)$$

$$v_{je} \geq u_{je} - u_{j,(e-1)} \quad \forall j, e \geq 2 \quad (13)$$

$$v_{je} = u_{je} \quad \forall j, e = 1 \quad (14)$$

### B. TEP-SNO-T1

The above sub-section presents a typical formulation for a traditional TEP problem. However, it fails to consider the flexibility in the transmission network since it assumes that all available transmission elements are connected in the grid. This is not necessary as treating transmission network as a dynamic network may promote more efficient solutions for delivering cheap power to load pockets. As load may have seasonal patterns, this sub-section improves the traditional TEP model by incorporating seasonal transmission network optimization. The first type of TEP-SNO only optimizes existing network without considering new lines. This change can be achieved by replacing (7)-(8) with (15)-(17). Variable $z_{kse}$ denotes the status of existing line $k$: 1 if it is connected to the grid; or 0 if it is disconnected from the network. Constraint (15) shows the branch capacity limit while (16)-(17) defines the power flow on existing branch $k$.

$$-z_{kse}RateA_k \leq p_{ktse} \leq z_{kse}RateA_k \quad \forall k,t,s,e \quad (15)$$

$$p_{ktse} - \frac{\theta_{n_{(k-)}tse} - \theta_{n_{(k+)}tse}}{x_k} \leq (1 - z_{kse})M \quad \forall k,t,s,e \quad (16)$$

$$p_{ktse} - \frac{\theta_{n_{(k-)}tse} - \theta_{n_{(k+)}tse}}{x_k} \geq -(1 - z_{kse})M \quad \forall k,t,s,e \quad (17)$$

### C. TEP-SNO-T2

TEP-SNO-T1 enhances the traditional TEP model by including the flexibility in the current network. This sub-section presents TEP-SNO-T2 that extends TEP-SNO-T1 by considering the flexibility provided by optimizing new lines. This means that a new line after its construction may be disconnected at a later point in time to improve the grid operational efficiency. This may sound counterintuitive since disconnection of a new line may discourage its initial purchase decision. However, it is possible that a new line would relieve network congestion for some scenarios while it worsens network congestion for some other scenarios. Thus, TEP-SNO-T2 can potentially provide a better solution than TEP-SNO-T1. To model the flexibility provided by new lines, (9)-(11) can be replaced by (18)-(21). Variable $z_{jse}$ denotes the status of new line $j$. Constraint (18) enforces capacity limits of candidate new lines; and (19)-(20) calculate the flow on new lines. Constraint (21) defines the relationship between $z_{jse}$ and $u_{je}$, which ensures that a new line cannot be connected to the grid if it is not available.

$$-z_{jse}RateA_j \leq p_{jtse} \leq z_{jse}RateA_j \quad \forall j,t,s,e \quad (18)$$

$$p_{jtse} - \frac{\theta_{n_{(j-)}tse} - \theta_{n_{(j+)}tse}}{x_j} \leq (1 - z_{jse})M \quad \forall k,t,s,e \quad (19)$$

$$p_{jtse} - \frac{\theta_{n_{(j-)}tse} - \theta_{n_{(j+)}tse}}{x_j} \geq -(1 - z_{jse})M \quad \forall k,t,s,e \quad (20)$$

$$z_{jse} \leq u_{je} \quad \forall j,s,e \quad (21)$$

### D. Various Planning Models

The traditional TEP model and two enhanced TEP-SNO models proposed in this paper are summarized in Table I. Existing transmission elements are treated as static assets in the planning horizon in the traditional TEP model; all new lines are also considered to be connected in the system after their construction is complete. The traditional TEP model consists of (1)-(14), which represents a mixed integer linear programming (MILP) problem. It is implemented to serve as a benchmark to demonstrate the cost saving benefits provided by incorporating SNO into the TEP model. As compared to the traditional TEP model, TEP-SNO-T1 includes the flexibility in existing network by replacing (7)-(8) with (15)-(17). TEP-SNO-T2 extends TEP-SNO-T1 by including the flexibility on new lines by replacing (9)-(11) with (18)-(21).

Table I. Summary of various TEP models

| | Constraints on existing lines | Constraints on new lines | Objective and other constraints |
|---|---|---|---|
| TEP | (7)-(8) | (9)-(14) | (1)-(6) |
| TEP-SNO-T1 | (15)-(17) | (9)-(14) | |
| TEP-SNO-T2 | (15)-(17) | (12)-(14), (18)-(21) | |

### E. Metrics:

It is very important to define the metrics for measuring the effectiveness of the proposed planning models. One key metric is to compare the cost of the proposed TEP-SNO models with the traditional TEP benchmark model. The absolute and relative cost reductions achieved with SNO are defined in (22) and (23) respectively.

$$TCR^{SNO} = TC^{TEPSNO} - TC^{TEP} \qquad (22)$$

$$\rho^{SNO} = TCR^{SNO}/TC^{TEP} \qquad (22)$$

where $TCR^{SNO}$, $TC^{TEPSNO}$, and $TC^{TEP}$ denote the cost reduction with SNO, total cost for TEP-SNO, and total cost for TEP respectively; $\rho^{SNO}$ denotes the percentage cost reduction with SNO.

## III. CASE STUDIES

The proposed traditional TEP and two enhanced TEP-SNO models are tested on a modified 24-bus test case that represents an area of the IEEE RTS-96 system [19]. This test case has 24 buses, 38 branches and 33 generators. The network topology, bus numbers and branch numbers for this test system are the same with the case used in [20]. The total system load in different intervals in typical days of different seasons is shown in Fig. 1. As shown in this figure, we can observe that summer is the peak season, followed by winter. The load profiles for spring season and fall season are very similar; both are lower than winter and summer. In this paper, the candidate new line set $J$ is a subset of existing line set $K$, which indicates that each candidate new line would have a parallel existing line. For convenience, the index of candidate new line set $J$ is consistent with existing line set $K$: if index $j$ has the same value with index $k$, then, new line $j$ connects to the same two buses that are connected by existing line $k$.

The computer that is used to conducted the simulations in this paper has the following specifications: Windows 10 Enterprise 64-bit, Intel(R) Core(TM) i7-8850H CPU 2.60GHz (12 CPUs), and 32GB memory. Gurobi 8.1.0 is used as the MILP solver to solve the transmission expansion planning problems; its termination condition is: 0.001% for *mipgap* or 3,000 seconds for *timelim*. AMPL is used to implement the model and serves as an interface between the model and the solver. In this paper, 3 epochs are considered in the planning horizon and the length of each epoch is 5 years.

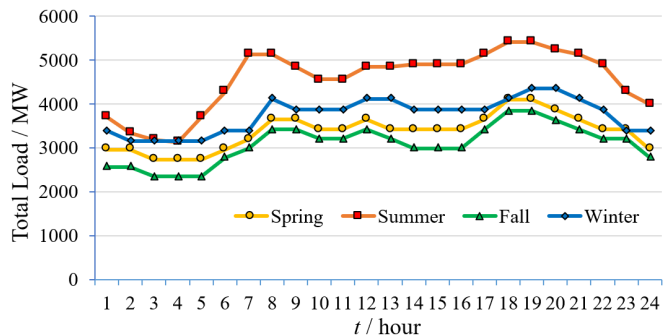

Fig. 1. Total load of the IEEE 24-bus test system in different hours in different seasons.

Table II presents the simulations results with various planning models. It shows that including SNO in the TEP model can significantly reduce the total planning cost. As compared to the traditional TEP model, the cost reduction with TEP-SNO-T1 is 465.0 million dollars or 4.69%; and the cost reduction with TEP-SNO-T2 is 491.2 million dollars or 4.95%. When compared with TEP-SNO-T1, TEP-SNO-T2 still achieves a lower cost by the amount of 26.2 million dollar. It is worth noting that the solutions obtained with the proposed TEP-SNO models correspond to a relative mipgap of 2.24% and 1.76% respectively while the TEP solution is optimum that has a mipgap of zero. This indicates that better solutions may be obtained with TEP-SNO if a higher time limit is used.

Table II. Simulation results with various planning models

|  | TEP | TEP-SNO-T1 | TEP-SNO-T2 |
|---|---|---|---|
| Total cost ($) | 9,921,190,000 | 9,456,150,000 | 9,429,990,000 |
| $TCR^{SNO}$ ($) | N/A | 465,040,000 | 491,200,000 |
| $\rho^{SNO}$ | N/A | 4.69% | 4.95% |
| Investment cost ($) | 1,367,480,000 | 683,981,000 | 683,981,000 |
| Generation cost ($) | 8,553,710,000 | 8,772,170,000 | 8,746,010,000 |
| Relative MipGap | 0.0000% | 2.2408% | 1.7646% |
| Best bound ($) | 9,921,190,000 | 9,244,260,000 | 9,263,590,000 |
| Solution time (s) | 8.66 | 3,000.31 | 3,000.48 |

Table III shows the transmission investment decisions with different planning models. The traditional TEP model suggests a plan of 6 new line construction: 4 new lines in epoch 1, and 2 new lines in epoch 2. The proposed two TEP-SNO models suggest the same planning decision: 3 new lines in epoch 1, and 1 new line in epoch 3. Thus, we can conclude that incorporating SNO into TEP can defer the construction of new lines and avoid unnecessary purchase of some other new lines. This explains why the investment cost for TEP-SNO is much lower than the traditional TEP model as shown in Table II.

Table III. Investment decisions with various planning models

|  | # of new lines | Epoch 1 | Epoch 2 | Epoch 3 |
|---|---|---|---|---|
| TEP | 6 | 13, 19, 22, 23 | 7, 27 | N/A |
| TEP-SNO-T1 | 4 | 19, 22, 23 | N/A | 11 |
| TEP-SNO-T2 | 4 | 19, 22, 23 | N/A | 11 |

N/A means that no new line is constructed in that epoch.

Table IV and Table V present the SNO results on existing lines with TEP-SNO-T1 and TEP-SNO-T2 respectively. For both TEP-SNO-T1 and TEP-SNO-T2, a total of 64 existing lines are disconnected for 4 seasons in 3 epochs in the planning horizon; note that the 64 transmission switching actions are very different for the two TEP-SNO models. In addition to switching existing lines off, the solution with TEP-SNO-T2 also switches new line 11 off in the fall season in epoch 3 that is after two seasons of its construction. This indicates that new line 11 can improve the network efficiency for Spring, Summer and Winter but not for Fall.

As shown in Table VI, the traditional TEP problem optimizes a model with 34,070 rows, 31,478 columns and 130,458 nonzeros that involves 31,394 continuous variables and 84 binary variables; the enhanced TEP-SNO-T1 optimizes a model with 65,846 rows, 31,912 columns and 225,786 nonzeros that involves 31394 continuous variables and 518 binary variables; the enhanced TEP-SNO-T2 optimizes a model with 66,014 rows, 32,080 columns and 226,122 nonzeros that involves 31,394 continuous variables and 686 binary variables. This statistics information indicates that the proposed two enhanced planning models are much more complex than the traditional TEP model, which is consistent with the solution time for these three planning models shown in Table II. Obviously, there is a tradeoff between solution time and solution quality. Though planning problems are not very solution time sensitive, a promising future work is the development of fast effective decomposition algorithm that can accelerate the solution process for solving the proposed enhanced TEP-SNO models.

Table IV. SNO results on existing lines with TEP-SNO-T1

| Season | Epoch 1 | Epoch 2 | Epoch 3 |
|---|---|---|---|
| Spring | 4, 31, 32, 33, 36, 37 | 15, 36, 38 | 1, 32, 33, 34, 35, 38 |
| Summer | 33, 34, 35, 37 | 6, 32, 35, 37 | 4, 6, 33, 34, 36 |
| Fall | 4, 31, 32, 33, 34, 35, 37 | 1, 15, 31, 32, 33, 36, 37 | 1, 11, 33, 34, 35, 37 |
| Winter | 3, 14, 34, 35 | 4, 20, 34, 36, 37, 38 | 1, 4, 20, 36, 37, 38 |

Table V. SNO results on existing lines with TEP-SNO-T2

| Season | Epoch 1 | Epoch 2 | Epoch 3 |
|---|---|---|---|
| Spring | 1, 15, 31, 32, 33, 34, 36, 37 | 1, 2, 32, 33, 38 | 1, 31, 32, 34, 35 |
| Summer | 6, 31, 34 | 4, 6, 33, 34, 35 | 4, 6, 33, 35, 36 |
| Fall | 9, 31, 33, 37 | 1, 15, 32, 33, 34, 35, 37, 38 | 1, 4, 32, 34, 35, 36, 38 |
| Winter | 3, 31, 34, 35 | 4, 32, 33, 34, 35 | 4, 32, 33, 34, 36 |

Table VI. Statistics of computational complexity for various planning models

|  | TEP | TEP-SNO-T1 | TEP-SNO-T2 |
|---|---|---|---|
| # of variables | 31,478 | 31,912 | 32,080 |
| # of linear variables | 31,394 | 31,394 | 31,394 |
| # of binary variables | 84 | 518 | 686 |
| # of constraints | 34,070 | 65,846 | 66,014 |
| # of equality constraints | 17,872 | 6,928 | 6,928 |
| # of inequality constraints | 16,198 | 58,918 | 59,086 |
| # of nonzeros in constraints | 130,458 | 225,786 | 226,122 |

IV. CONCLUSION

As load grows over time, current transmission network capacity would become a severe limiting factor for efficient and reliable grid operations. Thus, it is important to conduct TEP studies to make effective decisions on new transmission element investment. Traditional TEP model does not dynamically optimize the transmission network. However, seasonal load patterns may result in different optimal network configurations. Thus, this paper incorporates seasonal network optimization into the TEP model to utilize the flexibility in the transmission network. The proposed TEP-SNO-T1 model optimizes the status of existing branches while the proposed TEP-SNO-T2 model optimize the status of both existing branches and new branches. Numerical simulations show that both of the proposed two enhanced TEP-SNO models can improve system operational efficiency, defer investment of new transmission elements, and reduce the total cost. As compared to TEP-SNO-T1, TEP-SNO-T2 provides a better solution as it incorporates additional flexibility from optimizing the new line status.